# Numerical illustration of a recently proposed strongly polynomial-time algorithm for the general LP problem

Samuel Awoniyi
Department of Industrial and Manufacturing Engineering
FAMU-FSU College of Engineering
2525 Pottsdamer Street
Tallahassee, FL 32310
E-mail: awoniyi@eng.famu.fsu.edu
ORCID# 0000 0001 7102 6257

.

## Abstract

This article presents a numerical illustration of a recently proposed strongly polynomial-time algorithm for the general linear programming (LP) problem. This article is essentially the first half of an article that describes the proposed algorithm. Each iteration of the algorithm consists of two Gauss-Jordan pivoting operations. The algorithm is terminated after at most 2(k+n) iterations, where *k* is the number of constraints of the LP problem and *n* is the number of variables. Illustrative example LP problems described in this article include a Klee-Minty LP problem and an LP problem of Beale.

.

## 1. Introduction

In currently established literature on linear programming (LP), to find a strongly polynomial-time algorithm for the general LP problem is still an open problem [6,7,8,16]. Classical references on this topic include [2,3,9,10,11,13,15]. This article presents an illustration of a recently proposed strongly polynomial-time algorithm for solving the general LP problem [18]. The article [18] is quite long, and this illustration article is essentially a concise version of the first half of the article [18].

The proposed algorithm begins by utilizing LP duality theory to translate solving the general LP problem, having $k$ inequality constraints and $n$ variables, into solving a special system of equations in $R^{2(k+n)}$. Each iteration consists of two complementary Gauss-Jordan reduction pivoting in that system of equations, and does not require checking that problem variables are non-negative along the way as in the simplex algorithm. Each run of the algorithm stops after at most $2(k + n)$ iterations.

The rest of this article is organized as follows. Section 2 gives a problem statement wherein the general LP problem is translated into an equation-solving problem denoted as (Eq). This is followed in Section 3 by a statement of the proposed algorithm's iterations step-by-step. Section 4 presents the result of applying the proposed algorithm to some illustrative LP problem instances.

.

## 2. Problem statement

In this article, vectors are column vectors unless otherwise indicated. Vectors will be denoted by lower-case letters, and matrices by upper-case letters. Superscript $T$ will denote vector or matrix transpose as usual, and $I_{(.)}$ is reserved for identity matrix of dimension indicated by (.).

.

It is assumed that the general LP problem is given in Neumann symmetric form, (P) below:

.

$$\left\{\begin{array}{ll} \max & c^T x \\ \text{ST} & Ax \le b, \\ & x \ge 0 \end{array}\right\} \cdots\cdots\text{(P)}$$

.

where $c$ is $n$-vector, $A$ is $k$-by-$n$ matrix, $b$ is $k$-vector, and $x$ is $n$-vector of the problem's variables. From basic LP duality theory, solving (P) is equivalent to computing a $2(k+n)$-vector $z$ that solves the system of equations (Eq) stated below:

.

$$\left\{\begin{array}{ll} & Mz = q, \\ \text{ST} & z \ge 0 \ \& \\ & z_j z_{(m+n+j)} = 0, \ \text{for } j = 1, \ldots, m+n. \end{array}\right\} \cdots\cdots\text{(Eq)}$$

.

where

$$M = \begin{pmatrix} O & A & I_{(m)} & O \\ -A^T & O & O & I_{(n)} \\ -b^T & c^T & o^T & o^T \end{pmatrix} \text{ and } q = \begin{pmatrix} b \\ -c \\ o \end{pmatrix}$$

.

As an illustration of this problem formulation, let

.

$$c = \begin{pmatrix} -1 \\ 1 \end{pmatrix}, \ A = \begin{pmatrix} 1 & 1 \\ -1 & 0 \end{pmatrix}, \ b = \begin{pmatrix} 10 \\ -5 \end{pmatrix}$$

.

Then, in this instance, we have

.

$$M = \begin{bmatrix} 0 & 0 & 1 & 1 & 1 & 0 & 0 & 0 \\ 0 & 0 & -1 & 0 & 0 & 1 & 0 & 0 \\ -1 & 1 & 0 & 0 & 0 & 0 & 1 & 0 \\ -1 & 0 & 0 & 0 & 0 & 0 & 0 & 1 \\ 10 & -5 & 1 & -1 & 0 & 0 & 0 & 0 \end{bmatrix} \text{ and } q = \begin{bmatrix} 10 \\ -5 \\ 1 \\ -1 \\ 0 \end{bmatrix}$$

.

Problem (Eq) is an instance of primal-dual formulation of the general LP problem (see [6] for examples). A widely studied primal-dual formulation known as the Linear Complementarity Problem (LCP) (see [3,4], for example) is what one gets from (Eq) by not including the last equation of $Mz = q$. As it is well-known that the LCP is equivalent to the general LP problem in terms of solution existence, one might then surmise that the last row of matrix $M$ is redundant information. But Lemma 6.1 in [18] shows that the last row of $M$ is indeed quite indispensable.

.

# 3. Steps of the proposed algorithm for solving (Eq)

The proposed algorithm is a special pivoting method for solving (Eq), and the steps are as follows.

.

1 - **Initialize** the algorithm (as described below)
2 - **Stop if** stopping condition (described below) is met; **otherwise** go to step 3
3 - **Execute next iteration** (as described below), and thereafter go back to step 2 above

.

## 3.1 Initialization

Initialization consists of setting up an initial organizer/tableau for the algorithm's iterations. As notation henceforth, we let $[M\ q]$ denote the augmented matrix combining matrix $M$ and column vector $q$ ($M$ and $q$ as introduced in problem (Eq) in Section 2).

.

*Initialization operation:* Add the $(k + n + 1)$-th row of $[M\ q]$ to every other row of $[M\ q]$, in order to facilitate needed complementary pivoting on diagonal elements.

.

Let us temporarily denote the resultant matrix by $[\overline{M}\ \overline{q}]$. As numerical illustration, for our illustration example of Section 2 we have

.

$$[\overline{M}\ \overline{q}] = \begin{bmatrix} 10 & -5 & 2 & 0 & 1 & 0 & 0 & 0 & 10 \\ 10 & -5 & 0 & -1 & 0 & 1 & 0 & 0 & -5 \\ 9 & -4 & 1 & -1 & 0 & 0 & 1 & 0 & 1 \\ 9 & -5 & 1 & -1 & 0 & 0 & 0 & 1 & -1 \\ 10 & -5 & 1 & -1 & 0 & 0 & 0 & 0 & 0 \end{bmatrix}$$

.

From basic linear algebra, one can see that the system of equations and inequalities associated with $[\overline{M}\ \overline{q}]$ is equivalent to (Eq) in terms of solution existence. Accordingly, in the interest of notation tidiness, we will generally write $[M\ q]$ in place of $[\overline{M}\ \overline{q}]$ or any of its equivalent systems that we will obtain through elementary row operations on $[\overline{M}\ \overline{q}]$.

As a part of this initialization, define $\Pi$ as the set of indices of columns of $M$ that are *known* to be utilized in a solution basis matrix or *known* not to be utilized by any solution of (Eq), assuming that (Eq) has a solution. The index set $\Pi$ is initialized as the empty set $\varnothing$. Later, $\Pi$ will be updated continually, as described under "Executing next iteration" below.

In the remainder of this article, it may be advisable to keep in view the following general form of the matrix $[M\ q]$ :

| $m_{1,1}$ | $m_{1,2}$ | $\cdots$ | $m_{1,2(k+n)}$ | $q_1$ |
|---|---|---|---|---|
| $\vdots$ | $\vdots$ | $\ddots$ | $\vdots$ | $\vdots$ |
| $m_{k+n,1}$ | $m_{k+n,2}$ | $\cdots$ | $m_{k+n,2(k+n)}$ | $q_{k+n}$ |
| $m_{k+n+1,1}$ | $m_{k+n+1,2}$ | $\cdots$ | $m_{k+n+1,2(k+n)}$ | $q_{k+n+1}$ |

.

## 3.2 Stopping

There are two types of stopping - the case when a solution for (Eq) is found, and the case when there is evidence that (Eq) has no solutions.

.

*Case 1*: A solution of (Eq) is found
A solution of (Eq) is indicated in current $[M\ q]$ by having $q \geq 0$ along with $q_{k+n+1} = 0$.

.

*Case 2*: There is evidence that (Eq) has no solutions

A lack of solutions for (Eq) is indicated in current $[M\ q]$ by having $q_{k+n+1} > 0$ along with all other elements non-positive (that is, $\leq 0$) in row k+n+1 of $[M\ q]$ (that possibly after we multiply row $\mathrm{k+n+1}$ by $-1$ to effect $q_{k+n+1} > 0$). A lack of solutions is also indicated through a complementarity theorem contradiction that can occur in Step 1 or Step 4 of the description of MinorP pivoting and MajorP pivoting (under "Executing next iteration") below.

.

### 3.3 Executing next iteration

Towards describing "Executing next iteration", we define, in (i),..., (iv) below, several concepts that will feature in the remainder of this article.

.

`Definitions`

.

(i) A GJ pivoting in column $j$ of $[M\ q]$ is called *the complementary GJ pivoting in column j* if the pivoting position in $[M\ q]$ is (j,j) for $j \leq k+n$, or $(\mathrm{j-k-n,j})$ for $j > k+n$.

(ii) For $j \leq k+n$, column $j$ is the *complement column* for column k+n+j, and vice versa. Also, the index $j$ may be refered to as complement for the index $\mathrm{k+n+j}$, and vice versa.

(iii) Each iteration of the algorithm consists of two complementary GJ pivoting instances – a Minor Pivoting (abbreviated as *MinorP*) instance, when $q_{k+n+1} = 0$, followed by a Major Pivoting (abbreviated as *MajorP*) instance, when $q_{k+n+1} > 0$ (or, equivalently, $q_{k+n+1} < 0$). MinorP pivoting and MajorP pivoting instances will be described shortly in this section.

(iv) We will sometimes refer to a column of $M$ in $[M\ q]$, say $M^{(s)}$, as a *maximal column* if $m_{k+n+1,s} \geq m_{k+n+1,j}$ for $j = 1,\ldots,k+n$; that is, its $(\mathrm{k+n+1})$-th component (its element in the last row of $M$) is not smaller than that of any other column of $M$.

.

`Illustration of definitions`

.

The two concepts defined in (i) and (ii) are illustrated by the following table.

<table>
<tr><td>*</td><td></td><td></td><td></td><td>*</td><td></td><td></td><td></td><td>q<sub>1</sub></td></tr>
<tr><td></td><td>*</td><td></td><td></td><td></td><td>*</td><td></td><td></td><td>q<sub>2</sub></td></tr>
<tr><td></td><td></td><td>*</td><td></td><td></td><td></td><td>*</td><td></td><td>q<sub>3</sub></td></tr>
<tr><td></td><td></td><td></td><td>*</td><td></td><td></td><td></td><td>*</td><td>q<sub>4</sub></td></tr>
<tr><td colspan="8">last row (row 5) of <i>M</i></td><td>q<sub>5</sub></td></tr>
</table>

This table displays a frame of $[M\ q]$ for $k = n = 2$. An asterisk in the (i,j)-th position indicates the complementary pivoting position in column $j$, for $j = 1,...,8$. In each row, there are asterisks in two positions, and the columns of the two positions are complement columns for each other.

.

`A Cursory View on MinorP pivoting and MajorP pivoting`

.

Recall that each iteration consists of a MinorP pivoting instance followed by a MajorP pivoting instance. Each iteration utilizes MajorP to gradually assemble a basis matrix for a solution of (Eq), at the rate of about one column selection per MajorP instance. MinorP instances are intended to aid MajorP instances in that task, and can possibly make correct column selections as well.

Each column selection (by MajorP and MinorP instances) is from among columns of current $[M\ q]$ that have positive (k+n+1)-th component (that is, in the last row of $[M\ q]$). The column selections are carefully done (in Step 3 below) such that a kind of reversal of a previous MajorP column selection would promptly (in Step 4 below) yield a solution of (Eq), if (Eq) has any solutions.

A few details have been omitted in the foregoing cursory view of MinorP and MajorP. We will now fully describe "MinorP pivoting instance" and "MajorP pivoting instance". Each of the two types of pivoting instance is described in a four-step procedure.

.

MinorP pivoting instance (Here $q_{k+n+1} = 0$ and $q_i < 0$ some $i \leq k+n$)

If necessary, multiply the last row of $[M\, q]$ (that is, row k+n+1) by –1 to ensure that each negative component of $q$ (in $[M\, q]$) corresponds to a positive component of the last row of $M$. More precisely, multiply the last row of $[M\, q]$ (that is, row k+n+1) by –1 to ensure that $q_i(m_{k+n+1,i} + m_{k+n+1,k+n+i}) < 0$ for every $i \leq k+n$ having $q_i < 0$. The article [19] explains the feasibility of doing that.

*Step 1: Let $\mathcal{L}$ be the ordered list of column indices $j$ having $m_{k+n+1,j} > 0,$ in descending order of $m_{k+n+1,j} > 0.$ Let L denote $\mathcal{L}\backslash\Pi.$ The items in L are to be "picked up" one-at-a-time for processing. If $\mathrm{L} = \varnothing,$ then current iteration declares that (Eq) has no solutions by virtue of an LP strict complementarity theorem, as $\mathcal{L} \neq \varnothing$; otherwise, there are two cases to consider here.

.. Case 1.1: L has exactly 1 item in it – In this case, letting $\mathrm{L} = \{\widehat{j}\},$ current iteration (i) updates $\Pi$ by putting $\widehat{j}$ and its complement index into $\Pi,$ and (ii) performs desired MinorP pivoting in $\widehat{j}$-th column.

.. Case 1.2: L does not have exactly 1 item in it – In this case, with L having at least 2 elements in it, current iteration goes to Step 2.

*Step 2: Current iteration processes the next not-yet-processed item in list L. There are two cases to consider here.

.. Case 2.1: Some items are still available in list L to be processed – In this case, current iteration "picks up" the next item in L, labels the corresponding column as "current potential column selection", and then goes to Step 3.

..Case 2.2: There are no not-yet-processed items in L – In this case, current iteration goes to Step 4 where "finalizing operations" will be performed on the set L.

*Step 3: Current iteration processes the given "current potential column selection". There are two cases.

.. Case 3.1: The given "current potential column selection" is not the complement column for a previous MajorP "column selection" – In this case, the "current potential column selection" is labelled as "column selection" for this instance of MinorP, and current iteration thereafter performs desired MinorP pivoting in the column just labelled as "column selection".

.. Case 3.2: The given "current potential column selection" is the complement column for a previous MajorP "column selection" – In this case, the iteration goes back to Step 2.

*Step 4: First, do Task 1, Task 2 & Task 3 below, for each $\widehat{j} \in \mathrm{L}$ separately, until *either* a solution of (Eq) is obtained, *or* the current iteration is ended with $\widehat{j}$ put in the index set $\Pi$.

Task 1: As a probe, imagine/suppose that MinorP pivoting is to be performed in column $\widehat{j}$, followed by pertinent MajorP pivoting, say in column $w$.

Task 2: If that MajorP pivoting in $w$ would yield a solution of (Eq), then perform that supposed MinorP pivoting in $\widehat{j}$ and that supposed MajorP pivoting in $w,$ and then current iteration terminates the algorithm at this juncture. Otherwise, go to Task 3.

Task 3: In place of the pivoting sequence imagined in Task 1, consider doing a complementary GJ pivoting in column $w$ first, followed by doing a complementary pivoting in column $\widehat{j}$. If that reversal would indicate that $\widehat{j}$ should be put in the index set $\Pi,$ then do those two complementary GJ pivotings (that is, do complementary GJ pivoting in column $w$ first, followed by doing it in column $\widehat{j}$), put $\widehat{j}$ and its complement in $\Pi,$ and terminate the current iteration at this juncture.

But, if doing Tasks 1, 2, 3 for elements of L exhausts the set L without yielding a solution of (Eq) as in Task 2, and without putting $\widehat{j}$ in $\Pi$ and ending the current iteration as in Task 3, then the algorithm is to be terminated at this juncture, along with the declaration that (Eq) has no solutions (on account of a complementarity theorem contradiction).

.

MajorP pivoting instance (Here $q_{k+n+1} > 0$ or $q_{k+n+1} < 0$)

If $q_{k+n+1} < 0,$ then implicitly multiply row $k + n + 1$ of $[M\ q]$ by $-1$, to have $q_{k+n+1} > 0$.

*Step 1: Let $\mathcal{L}$ be the ordered list of column indices $j$ having $m_{k+n+1,j} > 0,$ in descending order of $m_{k+n+1,j} > 0.$ Let L denote $\mathcal{L}\backslash\Pi.$ The items in L are to be "picked up" one-at-a-time for processing. If $L = \varnothing,$ then current iteration declares that (Eq) has no solutions by virtue of an LP complementarity theorem, as $\mathcal{L} \neq \varnothing$; otherwise, there are two cases to consider here.

.. Case 1.1: L has exactly 1 item in it – In this case, letting $L = \{\hat{j}\},$ current iteration (i) updates $\Pi$ by putting $\hat{j}$ and its complement index into $\Pi,$ and (ii) performs desired MajorP pivoting in $\hat{j}$-th column.

.. Case 1.2: L does not have exactly 1 item in it – In this case, with L having at least 2 elements in it, current iteration goes to Step 2.

*Step 2: Current iteration processes the next not-yet-processed item in list L. There are two cases.

.. Case 2.1: Some items are still available in list L to be processed – In this case, current iteration "picks up" the next available item in L, labels the corresponding column as "current potential column selection" and then goes to Step 3.

.. Case 2.2: There are no not-yet-processed items in L – In this case, current iteration goes to Step 4 where "finalizing operations" will be performed on the set L.

*Step 3: Current iteration processes the given "current potential column selection". There are two cases.

.. Case 3.1: The given "current potential column selection" is not the complement column for a previous MajorP "column selection" – In this case, the "current potential column selection" is labelled as "column selection" for this instance of MajorP. Current iteration thereafter performs desired MajorP pivoting in the column just labelled as "column selection".

.. Case 3.2: The given "current potential column selection" is the complement column for a previous MajorP "column selection" – In this case, the iteration goes back to Step 2.

*Step 4: For each item $\hat{j}$ in L, current iteration performs MajorP pivoting in column $\hat{j}$ separately, and checks whether a solution of (Eq) has been obtained. If doing that yields a solution of (Eq), then current iteration terminates the algorithm with a solution of (Eq). Otherwise, if doing that exhausts L without yielding a solution of (Eq), then current iteration declares that (Eq) has no solutions, and terminates the algorithm at this juncture (on account of a complementarity theorem contradiction).

### 3.4 A flow chart and an illustration of the iterations

We display here a summarising flow chart and a numerical illustration of the algorithm's iteration described above.

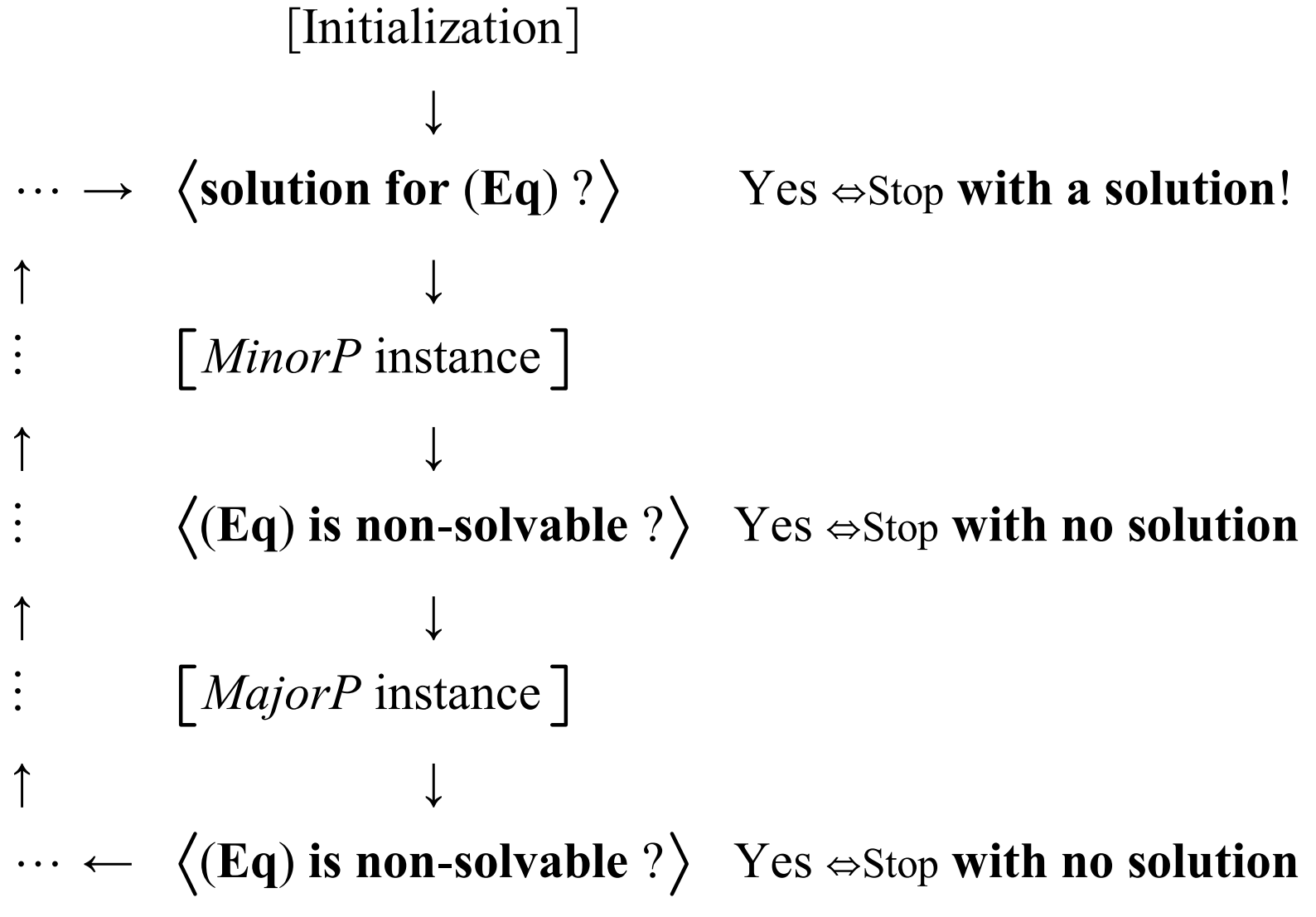

.

For a numerical illustration of the computations involved in the iterations, we will utilize the illustration LP example introduced in Section 2.

*Initialization*

$$[M\ q] = \begin{bmatrix} -10 & 5 & 0 & 2 & 1 & 0 & 0 & 0 & 10 \\ -10 & 5 & -2 & 1 & 0 & 1 & 0 & 0 & -5 \\ -11 & 6 & -1 & 1 & 0 & 0 & 1 & 0 & 1 \\ -11 & 5 & -1 & 1 & 0 & 0 & 0 & 1 & -1 \\ -10 & 5 & -1 & 1 & 0 & 0 & 0 & 0 & 0 \end{bmatrix}$$

Here, at the end of initialization, $\Pi$ is set to the empty set $\varnothing$.

.

*Iteration 1 MinorP*

L={2,4}. $\Pi = \varnothing$. Accordingly, the first instance of MinorP has its pivot at position (2,2), of current $[M\ q]$ instance to obtain the next $[M\ q]$ instance which we denote as $Z1$.

$$Z1 = \begin{bmatrix} 0 & 0 & 2 & 1 & 1 & -1 & 0 & 0 & 15 \\ -2 & 1 & -0.4 & 0.2 & 0 & 0.2 & 0 & 0 & -1 \\ 1 & 0 & 1.4 & -0.2 & 0 & -1.2 & 1 & 0 & 7 \\ -1 & 0 & 1 & 0 & 0 & -1 & 0 & 1 & 4 \\ 0 & 0 & 1 & 0 & 0 & -1 & 0 & 0 & 5 \end{bmatrix}$$

*Iteration 1 MajorP*

L={3}. On acccount of L being a singleton, $\Pi = \{3, 7\}$. Accordingly, the first instance of MajorP has its pivot at position (3,3) of current $[M\ q]$ instance, $Z1$, to obtain the next $[M\ q]$ instance which we denote as $P1$.

$$
\text{P1} = \begin{bmatrix} -1.4286 & 0 & 0 & 1.2857 & 1 & 0.7143 & 0 - 1.4286 & 0 & 5 \\ -1.7143 & 1 & 0 & 0.1429 & 0 & -0.1429 & 0.2857 & 0 & 1 \\ 0.7143 & 0 & 1 & -0.1429 & 0 & -0.8571 & 0.7143 & 0 & 5 \\ -1.7143 & 0 & 0 & 0.1429 & 0 & -0.1429 & -0.7143 & 1 & -1 \\ -0.7143 & 0 & 0 & 0.1429 & 0 & -0.1429 & -0.7143 & 0 & 0 \end{bmatrix}
$$

*Iteration 2 MinorP*

L={4}. On acccount of L being a singleton, $\Pi = \{3,7,4,8\}$. Accordingly, the second instance of MinorP has its pivot at position (4,4) of current $[M\ q]$ instance, $P1$, to obtain the next $[M\ q]$ instance which we denote as $Z2$.

$$
\text{Z2} = \begin{bmatrix} 14 & 0 & 0 & 0 & 1 & 2 & 5 & -9 & 14 \\ 0 & 1 & 0 & 0 & 0 & 0 & 1 & -1 & 2 \\ -1 & 0 & 1 & 0 & 0 & 0 & -1 & 1 & 4 \\ -12 & 0 & 0 & 1 & 0 & -1 & -5 & 7 & -7 \\ 1 & 0 & 0 & 0 & 0 & 0 & 0 & -1 & 1 \end{bmatrix}
$$

*Iteration 2 MajorP*

L={1}. On acccount of L being a singleton, $\Pi = \{3,7,4,8,1,5\}$. Accordingly, the second instance of MajorP has its pivot at position (1,1) of current $[M\ q]$, $Z2$, to obtain the next $[M\ q]$ instance which we denote as $P2$.

$$
\text{P2} = \begin{bmatrix} 1 & 0 & 0 & 0 & 0.0714 & 0.1429 & 0.3571 & -0.6429 & 1 \\ 0 & 1 & 0 & 0 & 0 & 0 & 1 & -1 & 2 \\ 0 & 0 & 1 & 0 & 0.0714 & -0.8571 & 0.3571 & 0.3571 & 5 \\ 0 & 0 & 0 & 1 & 0.8571 & 0.7143 & -0.7143 & -0.7143 & 5 \\ 0 & 0 & 0 & 0 & -0.0714 & -0.1429 & -0.3571 & -0.3571 & 0 \end{bmatrix}
$$

The algorithm stops here with a solution of (Eq), because $q \geq 0$ & $q_5 = 0$. Thus, the illustrative example problem introduced in Section 2 has been solved in two iterations of the algorithm; it does not require Step 4 of the procedures stated in Section 4.3.

.

## 3.5 Further clarification on the iterations

(i) In Step 1 of the 4-step procedure that describes MinorP pivoting, the ordering of index set L does not have to be "descending order" of $m_{k+n+1,j} > 0$. For example, it could be in "ascending order" of $m_{k+n+1,j} > 0$. Here is one important source of variation in the algorithm, and that presents some opportunity for utilizing special structures in special classes of LP problem classes. However, it is advisable to use the same ordering of L throughout all iterations of each run of the algorithm. In comparison to that, in Step 1 of the 4-step procedure that describes MajorP pivoting, the ordering of index set L has to be a descending order, on account of Lemma 6.1 which will be stated presently.

(ii) As a closing part of "column selection" in Case 3.1 of "MajorP pivoting instance" above, one can assume, without loss of generality, that the "column selection" (that is, the selected column) has maximal $(k + n + 1)$-th component, among all the columns indicated by L. That is without loss of generality, because the selected column may be multiplied by a positive number without altering the result of performing pertinent complementary GJ pivoting on that column.

(iii) One may describe our algorithm as consisting of an "initialization" procedure, followed by "MinorP → MajorP → MinorP cycles" that terminates with *either* a MajorP pivoting instance, together with a solution of (Eq), *or else* a MajorP or a MinorP pivoting instance indicating that (Eq) has no solutions. The MinorP →MajorP → MinorP cycles are a consequence of applying pairs of complementary GJ pivotings to skew-symmetric matrices, as explained in the article [19].

.

# 4. A report on some illustrative LP problems

We present in this Section a brief report on how our algorithm performed on some illustrative LP problems that include an instance of Klee-Minty LP problem and a Beale LP problem.

In our description of the illustrative example problems, we utilize what we call "Z-P Records Table". A Z-P Records Table displays the indices of columns of $[M\ q]$ that are selected for MinorP pivoting (that is, when $q_{k+n+1} = 0$, Z for zero) and MajorP pivoting (that is, when $q_{k+n+1} > 0$, P for positive). For the numerical example given under "An illustration of the iterations", the Z-P Records Table is:

| itn | Z | P |
|---|---|---|
| 1 | 4 | 1* |
| 2 | 2* | 3* |

We display results for two types of ordering of (the index set) $L$ for MinorP pivoting – Type 1 refers to the ascending ordering, and Type 2 refers to the descending ordering. In each example's Z-P Records Tables, column indices that belong to $\Pi$ are "asterisked" (as in the iteration illustration example just stated above).

Each example LP problem is assumed to be of the Neumann symmetric form

$$\begin{aligned} \text{maximize} \quad & f^T x \\ \text{subject to:} \quad & Ax \le b \\ & x \ge 0 \end{aligned}$$

with corresponding data given in the form

$$\left(\begin{array}{c|c} f^T & \\ \hline A & b \end{array}\right)$$

.

### Example 1: An instructive LP problem

Corresponding "Type 1 Z-P Records" (below) illustrates Step 4 of "MajorP pivoting instance" (described above) in iteration 5 – all candidates for a MajorP column selection (columns 4, 7 & 9) are the complements of previous MajorP column selections (columns 11, 14 & 2 respectively).

Similarly, "Type 2 Z-P Records" table illustrates Step 4 of "MinorP pivoting instance" in iteration 5 – all candidates for a MinorP column selection (columns 5 & 8) are the complements of previous MajorP column selections (columns 12 & 1 respectively).

In both case, the algorithm terminates there with a solution of (Eq), just as stated in corresponding Step 4.

$$\left(\begin{array}{cccc|c} 2 & 7 & 6 & 4 & \\ \hline 1 & 1 & 0.83 & 0.5 & 65 \\ 1.2 & 1 & 1 & 1.2 & 96 \\ 0.5 & 0.7 & 1.2 & 0.4 & 80 \end{array}\right);$$

with primal LP solution

$$x = (0, 5.1601, 53.2015, 31.3653)^T$$

and dual solution

$$y = (6.2147, 0.7062, 0.1130)^T$$

.

| Type 1 Z-P Records | | |
|---|---|---|
| itn | Z | P |
| 1 | 4 | 2 |
| 2 | 7 | 11 * |
| 3 | 6 | 3 |
| 4 | 5 * | 14 |
| 5 | 1 | 7 |

. .

| Type 2 Z-P Records | | |
|---|---|---|
| itn | Z | P |
| 1 | 5 | 1 |
| 2 | 7 | 2 |
| 3 | 6 * | 12 |
| 4 | 3 * | 4 |
| 5 | 5 | 11 |

Example 2: An unbounded LP problem

The algorithm indicates that (Eq) has no solutions, after a MinorP pivoting, just as expected,.

$$\left(\begin{array}{ccc|c} 1 & 2 & 1.5 & \\ \hline 1 & 0 & 0 & 5 \\ 0 & 1 & 0 & 6 \\ -1 & -1 & -1 & -10 \end{array}\right); \quad \begin{array}{c} \text{with primal LP solution} \\ \text{n. a.} \\ \text{and dual solution} \\ \text{n. a.} \end{array}$$

.

| Type 1 Z-P Records | | |
|---|---|---|
| itn | Z | P |
| 1 | 4 | 1 * |
| 2 | 6 | – |

. .

| Type 2 Z-P Records | | |
|---|---|---|
| itn | Z | P |
| 1 | 3 | 4 |
| 2 | 1 | 5 |
| 3 | 9 | 2 * |
| 4 | 6 | – |

.

Example 3: An instance of Beale LP problem

This is a classical LP problem that features a very degenerate primal solution. Our algorithm automatically rests on an LP strict complementarity theorem, so that it is not hampered by degenerate primal solution or dual solution. Examples 4, 6 & 7 illustrate that fact as well.

$$\left(\begin{array}{cccc|c} 0.75 & -150 & 0.02 & -6 & \\ \hline 0.25 & -60 & -0.04 & 9 & 0 \\ 0.50 & -90 & -0.02 & 3 & 0 \\ 0.00 & 0 & 1.00 & 0 & 1 \end{array}\right); \quad \begin{array}{c} \text{with primal LP solution} \\ x = (0.04, 0, 1, 0)^T \\ \text{and dual solution} \\ y = (0, 1.5, 0.05)^T \end{array}$$

.

| Type 1 Z-P Records | | |
|---|---|---|
| itn | Z | P |
| 1 | 6 | 3 ∗ |
| 2 | 4 ∗ | 2 |

. .

| Type 2 Z-P Records | | |
|---|---|---|
| itn | Z | P |
| 1 | 4 | 2 |
| 2 | 6 ∗ | 3 ∗ |

Example 4: An example with highly degenerate solutions for primal and dual

$$\left(\begin{array}{cccc|c} 2 & 1 & -1 & 1 & \\ \hline 1 & 1 & 1 & 1 & 12 \\ -1 & 0 & 1 & -1 & -8 \\ 0 & 2 & 0 & -1 & 6 \end{array}\right);$$

with primal LP solution $x = (12,0,0,0)^T$ and dual solution $y = (2,0,0)^T$

.

| Type 1 Z-P Records | | |
|---|---|---|
| itn | Z | P |
| 1 | 5 | 3 |
| 2 | 7 | 1* |
| 3 | 4 | 14 |
| 4 | 10* | 12* |

. .

| Type 2 Z-P Records | | |
|---|---|---|
| itn | Z | P |
| 1 | 2 | 4 |
| 2 | 9 | 1 |

Example 5: An instructive problem from p.57 of Dantzig's book [5]

The Type 2 MinorP pivotings of this example illustrate the usefulness of comparing the signs of components of $q$ and corresponding components of the last row of $[M\ q]$ whenever $q_{k+n+1} = 0$ (that is, prior to every MinorP pivoting).

.

$$\left(\begin{array}{ccccc|c} -2 & 1 & -3 & -7 & 5 & \\ \hline 1 & 2 & 1 & 1 & 6 & 10 \\ -2 & -3 & -4 & -1 & -2 & -4 \\ 3 & 2 & 0 & 3 & 1 & 8 \end{array}\right);$$

with primal LP solution $x = (0, 0.2857, 0, 0, 1.5714)^T$ and dual solution $y = (0.9286, 0.2857, 0)^T$

.

| Type 1 Z-P Records | | |
|---|---|---|
| itn | Z | P |
| 1 | 5 | 1 |
| 2 | 3 | 8 |
| 3 | 11 * | 2 |

. .

| Type 2 Z-P Records | | |
|---|---|---|
| itn | Z | P |
| 1 | 8 | 1 |
| 2 | 2 * | 6 |
| 3 | 5 * | 3 |
| 4 | 14 * | 7 |
| 5 | 4 | 15 * |
| 6 | 11 * | 12 |

```
Example 6: Another instructive example
```

This LP problem illustrates the iteration count 2(k+n) instead of k+n, as mentioned in the computational complexity lemma, Lemma 7.1 in [18].

$$\left(\begin{array}{cccc|c} 2 & -1 & 4 & 3 & \\ \hline 1 & 2 & 1 & 2 & 20 \\ 3 & 1 & -1 & 2 & 18 \\ 1 & 1 & -1 & -3 & 21 \end{array}\right);$$

with primal LP solution

$$x = (0,0,20,0)^T$$

and dual solution

$$y = (4,0,0)^T$$

.

| Type 1 Z-P Records | | |
|---|---|---|
| itn | Z | P |
| 1 | 4 | 2 |
| 2 | 7 | 1 |
| 3 | 11 | 6 |
| 4 | 4 | 3 |
| 5 | 14 | 5 |
| 6 | 11 | 10 * |
| 7 | 7 | 12 * |
| 8 | 4 | 14 * |
| 9 | 9 * | 11 |

. .

| Type 2 Z-P Records | | |
|---|---|---|
| itn | Z | P |
| 1 | 6 | 1 |

.

```
Example 7: An instance of Klee-Minty LP problem (with n=3)
```

This classical LP problem has a highly degenerate solution. Even though problem data is highly ill-conditioned in this example, our algorithm converges as predicted.

$$\left(\begin{array}{ccc|c} 100 & 10 & 1 & \\ \hline 1 & 0 & 0 & 1 \\ 20 & 1 & 0 & 100 \\ 200 & 20 & 1 & 10000 \end{array}\right);$$

with primal LP solution
$x = (0, 0, 10000)^T$
and dual solution
$y = (0, 0, 1)^T$

.

| Type 1 Z-P Records | | |
|---|---|---|
| itn | Z | P |
| 1 | 6 | 3 |

. .

| Type 2 (modified) Z-P Records | | |
|---|---|---|
| itn | Z | P |
| 1 | 4 | 3 |
| 2 | 2 | 5 |
| 3 | 10 | 6 |
| 4 | 4 | 1 |
| 5 | 8 ∗ | 11 ∗ |
| 6 | 7 ∗ | 10 |

.

The $n$-variable instance of Klee-Minty LP problem is solved by the Type 1 variant of our algorithm in exactly one iteration, with chosen pivot columns of $[M\ q]$ then being $(2n)$-th column by first MinorP instance, and $n$-th column by first MajorP instance.

.

# 5. Remarks on some further work

In the illustrative examples described above, the algorithm performed as expected in terms of number of iteration. However, one unexpected feature that became apparent is that ill-conditioned problem data, as in Examples 4, 6 &7, are adequately treated by implementing both Type 1 and Type 2 MinorP pivoting. Accordingly, one direction for further work is to do computational work to establish the efficacy of using variants of Type1 MinorP pivoting and Type 2 MinorP pivoting to address ill-conditioned data.

As the algorithm ordinarily generates a good deal of data, as it solves both primal and dual problems together, it should be of interest to investigate how such data could be helpful in solving classes of practical Mixed Integer Linear Programming (MILP) problems.

A compact implementation of the algorithm, similar to [20], could be a useful computational project as well.

.